\documentstyle[11pt]{amsart} \numberwithin{equation}{section}
\input{amssym.def}
\input{amssym.tex}

\begin{document}

\title[Lie bracket derivation-derivations]{Lie bracket derivation-derivations in Banach algebras}

\author[C. Park]{Choonkil Park}
\address{Choonkil Park \newline \indent Research Institute for Natural Sciences, Hanyang University, Seoul 04763,   Korea}
\email{baak@@hanyang.ac.kr}\vskip 2mm

\begin{abstract}
In this paper, we introduce and solve the following additive-additive  $(s,t)$-functional inequality
\begin{eqnarray}\label{0.1}  && \|g\left(x+y\right) -g(x) -g(y)\| +\| h(x+y) + h(x-y) -2 h(x) \|   \\ && \le \left\|s\left( 2 g\left(\frac{x+y}{2}\right)-g(x)-g(y)\right)\right\|+ \left\|t \left(  2h\left(\frac{x+y}{2}\right)+ 2h \left(\frac{x-y}{2}\right)- 2h (x)\right) \right\|   ,  \nonumber
\end{eqnarray}
where $s$ and $t$ are fixed nonzero complex numbers with $|s| <1$ and $ |t| <1$.
Using  the direct method and the fixed point method, we prove the Hyers-Ulam stability of Lie bracket derivation-derivations in complex Banach algebras, associated to the additive-additive $(s,t)$-functional inequality (\ref{0.1})  and the following functional inequality 
\begin{eqnarray} \label{0.2}\|  [g, h](xy)-[g,h](x) y- x [g,h](y) \| +\| h(xy) - h(x) y -x h(y) \| \le   \varphi(x,y).
\end{eqnarray}
\end{abstract}



\maketitle

\baselineskip=15pt

\theoremstyle{definition}
  \newtheorem{df}{Definition}[section]
  \newtheorem{rk}[df]{Remark}
\theoremstyle{plain}
  \newtheorem{lem}[df]{Lemma}
  \newtheorem{thm}[df]{Theorem}
  \newtheorem{pro}[df]{Proposition}
  \newtheorem{cor}[df]{Corollary}

\setcounter{section}{0}

\section{Introduction and preliminaries}

Let $A$ be a complex Banach algebra and  $Der(A)$ be the set of $\mathbb{C}$-linear bounded derivations on $A$.  For $\delta_1, \delta_2 \in  Der(A)$,  \begin{eqnarray*} \delta_1 \circ \delta_2 (ab) &=&\delta_1 \circ \delta_2 (a) b + \delta_2(a)\delta_1 (b)+  \delta_1(a) \delta_2(b) + a \delta_1 \circ \delta_2 (b) ,\\
\delta_2 \circ \delta_1 (ab) &=&\delta_2 \circ \delta_1 (a) b + \delta_1(a)\delta_2 (b) + \delta_2(a) \delta_1(b) + a \delta_2 \circ \delta_1 (b) 
\end{eqnarray*}
for all $a, b\in A$. Let $[\delta_1, \delta_2] = \delta_1 \circ \delta_2 - \delta_2 \circ \delta_1$. Then 
\begin{eqnarray*} [\delta_1 , \delta_2]  (ab) = [\delta_1 , \delta_2] (a) b +  a [\delta_1 , \delta_2] (b) 
\end{eqnarray*} for all $a,b \in A$. Since $[\delta_1 , \delta_2] : A \to A$ is $\mathbb{C}$-linear, $[\delta_1 , \delta_2] \in Der(A)$ for all $\delta_1, \delta_2 \in Der(A)$.
Thus $Der(A)$ is a Lie algebra with Lie bracket $[\delta_1, \delta_2]$, since $\delta_1 + \delta_2$ and $\alpha \delta_1$ are $\mathbb{C}$-linear derivations on $A$ for all $\delta_1, \delta_2 \in Der(A)$ and all  $\alpha \in \mathbb{C}$. One can easily show that $Der(A)$ is a Banach space, since $A$ is complete.

In this paper, we introduce and investigate  Lie bracket derivation-derivations  in  complex Banach algebras.

\begin{df} Let $A$ be a complex Banach algebra and $D, H: A\to A$ be $\mathbb{C}$-linear mappings. Let $[D,H](a) = D(H(a)) - H(D(a))$ for all $a\in A$. A $\mathbb{C}$-linear mapping $[D,H]:  A \to A$
is called a {\it  Lie bracket derivation-derivation} in $A$ if $[D, H]$ and $H$ (or $D$) are derivations in $A$, i.e., 
\begin{eqnarray*}
[D,H] (ab) & = & [D,H](a) b + a [D,H](b),\\ H(ab) &= &H(a) b + a H(b)
\end{eqnarray*}
for all $a,b \in A$. 
\end{df}

Since  $[\delta_1 , \delta_2] \in Der(A)$ for  $\delta_1, \delta_2 \in Der(A)$, $[\delta_1 , \delta_2]$ is a Lie bracket derivation-derivation.

The stability problem of functional equations  originated from a
question of Ulam \cite{ul60}  concerning the stability of group
homomorphisms.
 Hyers \cite{hy41} gave a first affirmative partial answer to the question of Ulam for Banach spaces.
Hyers' Theorem  was generalized by Aoki \cite{a} for additive
mappings and by  Rassias \cite{ra78} for linear mappings by
considering an unbounded Cauchy difference.  A
generalization of the  Rassias theorem was obtained by G\u
avruta \cite{ga94} by replacing the unbounded Cauchy difference by a
general control function in the spirit of  Rassias' approach.
  Park \cite{ppp, ppp1, par} defined additive $\rho$-functional inequalities and proved the Hyers-Ulam stability of the additive $\rho$-functional inequalities in Banach spaces and non-Archimedean Banach spaces.
 The
 stability problems of various functional equations and functional inequalities have been extensively investigated
 by a number of authors (see \cite{el, ega,   n15, sl, w}).  

We recall a fundamental result in fixed point theory.

\begin{thm} {\rm \cite{cr1, dm}}\label{thm1.1} Let $(X, d)$ be a complete generalized metric space and let $J: X \rightarrow X$ be a strictly contractive mapping
with Lipschitz constant $\alpha <1$. Then for each given element $x\in X$, either $$d(J^n x, J^{n+1} x) = \infty$$
for all nonnegative integers $n$ or there exists a positive integer $n_0$ such that

$(1)$ $d(J^n x, J^{n+1}x) <\infty, \qquad \forall n\ge n_0$;

$(2)$ the sequence $\{J^n x\}$ converges to a fixed point $y^*$ of $J$;

$(3)$ $y^*$ is the unique fixed point of $J$ in the set $Y = \{y\in X \mid d(J^{n_0} x, y) <\infty\}$;

$(4)$ $d(y, y^*) \le \frac{1}{1-\alpha} d(y, Jy)$ for all $y \in Y$.
\end{thm}

In 1996,  Isac and  Rassias \cite{ir} were the first to
provide applications of stability theory of functional equations for
the proof of new fixed point theorems with applications.
By using fixed point methods, the stability problems of several functional
equations have been extensively investigated by a number of
authors (see \cite{cr2, cr3, e17, pa17, r}).

In this paper, we solve the additive-additive $(s,t)$-functional inequality {\rm (\ref{0.1})}. Furthermore, we investigate Lie bracket derivation-derivations in complex Banach algebras associated to  the additive-additive $(s,t)$-functional inequality {\rm (\ref{0.1})} and the functional inequality (\ref{0.2}) by using the direct method and by the fixed point method.

Throughout this paper, assume that $A$ is a complex Banach algebra and that  $s$ and $t$ are fixed nonzero complex numbers with $|s|<1$ and $ |t| <1$.

\section{Stability of additive-additive $(s,t)$-functional inequality {\rm (\ref{0.1})}: a direct method}

In this section, we solve and  investigate the additive-additive  $(s,t)$-functional inequality (\ref{0.1}) in complex  Banach algebras.

\begin{lem}\label{lm2.1}
If  mappings  $g,h  : A \rightarrow A$ satisfy $g(0)=h(0)=0$ and 
\begin{eqnarray}\label{2.1} && \|g\left(x+y\right) -g(x) -g(y)\| +\| h(x+y) + h(x-y) -2 h(x) \|   \\ && \le \left\|s\left( 2 g\left(\frac{x+y}{2}\right)-g(x)-g(y)\right)\right\|+ \left\|t \left(  2h\left(\frac{x+y}{2}\right)+ 2h\left(\frac{x-y}{2}\right)- 2h(x)\right) \right\|   \nonumber
\end{eqnarray} for all $x, y \in A$, then the mappings   $g, h : A \rightarrow A$ are additive.
\end{lem}

\begin{pf}
Letting $y=x$ in {\rm (\ref{2.1})}, we get
\begin{eqnarray*}\|g(2x)- 2g(x)\| + \|h(2x)- 2h(x)\| \le 0
\end{eqnarray*} for all $x \in A$. 
So $g(2x) = 2g(x)$ and $h(2x)= 2h(x)$ for all $x\in A$.
It follows from (\ref{2.1}) that 
\begin{eqnarray*}   && \left\|g\left(x+y\right) -g(x) -g(y)\right\| + \left\| h\left(x+y\right) +h(x-y) -2 h(x) \right\|   \\ && \le \left\| s\left(  g\left(x+y\right)-g(x)-g(y)\right)\right\|+ \|t (  h(x+y)-h(x)-h(y))\| 
\end{eqnarray*} for all $x, y \in A$. 
Thus $g\left(x+y\right) -g(x) -g(y)=0$ and  $h\left(x+y\right) + h(x-y) - 2h(x) =0$ for all $x\in A$, since $|s|<1$ and $|t|<1$.
So the mappings  $g,h: A \to A $ are  additive.
\end{pf}

\begin{lem} \label{lm2.2} \cite[Theorem 2.1]{pa06}
Let  $f: A \to A$  be an additive   mapping such that $$f(\lambda a ) =  \lambda  f(a) $$  for all $\lambda \in { {\mathbb{T}}^1 } : = \{ \xi \in {\mathbb{C}} : | \xi | =1\}$ and all $a\in A$.  Then the mapping $f: A \to A$ is $\mathbb{C}$-linear.
\end{lem}

Using the direct method,  we prove the Hyers-Ulam stability of Lie bracket  derivation-derivations  in complex Banach algebras associated to  the additive-additive $(s,t)$-functional
inequality (\ref{2.1}) and the functional inequality (\ref{0.2}).

\begin{thm}\label{thm2.3}
Let $\varphi: A^2 \to [0,\infty)$ be a function such that 
\begin{eqnarray}\label{2.2}
 \sum_{j=1}^{\infty} 4^j \varphi\left(\frac{x}{2^j}, \frac{y}{2^j}\right)  <  \infty
\end{eqnarray}
for all $x, y\in A$. Let $g, h :
A \rightarrow A$ be   mappings  satisfying $g(0)=h(0)= 0$ and 
\begin{eqnarray}\label{2.3}
 &&\left\|g\left(\lambda (x+y)\right) -\lambda g(x) -\lambda g(y)\right\| + \left\| h\left(\lambda (x+y)\right)+h(\lambda  (x-y)) - 2\lambda h(x) \right\|  \nonumber  \\ &&\qquad  \le \left\| s\left( 2 g\left(\lambda \frac{x+y}{2}\right)-\lambda g(x)- \lambda g(y)\right)\right\|\\ && \qquad + \left\|t \left(  2h\left(\lambda \frac{x+y}{2}\right)+ 2 h\left(\lambda \frac{x-y}{2}\right)-2\lambda h(x)\right)\right\| +\varphi(x,y)    \nonumber
\end{eqnarray} 
for all $\lambda \in {{\mathbb{T}}^1} $ and all $x, y \in A$.  
If the mappings  $g, h: A\to A$ satisfy 
\begin{eqnarray} \label{2.4}
\|[g,h](xy) - [g,h](x) y - x [g,h](y) \| + \|h(xy) - h(x)y-xh(y)\| \le \varphi(x,y)
\end{eqnarray}
 for all $x,y\in A$, then there exist a unique  $\mathbb{C}$-linear  $D : A \rightarrow A$ and a unique  derivation  $H: A\to A$ such that $[D,H] : A\to A$ is a  derivation and 
\begin{eqnarray}\label{2.5}
\|g(x)-D(x)\| + \|h(x)-H(x)\| \le \sum_{j=1}^{\infty}  2^{j-1} \varphi\left(\frac{x}{2^j}, \frac{y}{2^j}\right) 
\end{eqnarray}
for all $x\in A$.
\end{thm}

\begin{pf}
Letting $\lambda =1$ and $y=x$  in  (\ref{2.3}), we get 
\begin{eqnarray}\label{2.6}
\|g(2x)-2g(x)\| + \|h(2x)-2h(x)\|\le \varphi(x,x)
\end{eqnarray}
and so 
\begin{eqnarray*}
\left\|g(x)-2g\left(\frac{x}{2}\right)\right\| + \left\|h(x)-2h\left(\frac{x}{2}\right)\right\|\le \varphi\left(\frac{x}{2},\frac{x}{2}\right)
\end{eqnarray*}
for all $x\in A$.  Thus  
\begin{eqnarray}\label{2.7}
&& \left\|2^{l} g\left(\frac{x}{2^l}\right) - 2^{m} g\left(\frac{x}{2^{m}}\right)\right\| + \left\|2^{l} h\left(\frac{x}{2^l}\right) - 2^{m} h\left(\frac{x}{2^{m}}\right)\right\|    \\ & &  \qquad \le 
 \sum_{j=l}^{m-1}\left\| 2^j g\left(\frac{x}{2^j}\right) - 2^{j+1} g\left(\frac{x}{2^{j+1}}\right) \right\| + \sum_{j=l}^{m-1}\left\| 2^j h\left(\frac{x}{2^j}\right) - 2^{j+1} h\left(\frac{x}{2^{j+1}}\right) \right\| \nonumber \\ && \qquad 
\le  \nonumber \sum_{j=l+1}^{m} 2^{j-1}  \varphi\left(\frac{x}{2^{j}},  \frac{x}{2^{j}}\right) \
\end{eqnarray}
 for all nonnegative
integers $m$ and $l$ with $m>l$ and all $x \in A$. It follows from
{\rm (\ref{2.7})} that the sequences
$\{2^{k} g(\frac{x}{2^k})\}$ and $\{2^{k} h(\frac{x}{2^k})\}$   are  Cauchy  for all $x \in A$.
Since $Y$ is a  Banach space, the sequences $\{2^{k} g(\frac{x}{2^k})\}$ and $\{2^{k} h(\frac{x}{2^k})\}$   converge.
So one can define
the mappings $D, H : A \rightarrow A$ by
$$D(x) : = \lim_{k\to \infty} 2^{k} g\left(\frac{x}{2^k}\right),  \qquad H(x) : = \lim_{k\to \infty} 2^{k} h\left(\frac{x}{2^k}\right)  $$
for all $x \in A$. Moreover, letting $l =0$ and passing to the limit $m \to
\infty$ in {\rm (\ref{2.7})}, we get {\rm (\ref{2.5})}.

It folllows from (\ref{2.3}) that \begin{eqnarray*}&& \|D\left(\lambda(x+y)\right) -\lambda D(x) - \lambda D(y)\| +\| H(\lambda(x+y)) + H(\lambda(x-y) )-2 \lambda H(x) \|   \\ && = \lim_{n\to\infty} 2^n \left\|g\left(\lambda \frac{x+y}{2^n}\right) -\lambda g\left(\frac{x}{2^n}\right) -\lambda g\left(\frac{y}{2^n}\right)\right\| \\ && + \lim_{n\to\infty} 2^n \left\| h\left(\lambda \frac{x+y}{2^n}\right) + h\left(\lambda \frac{x-y}{2^n}\right) -2 \lambda h\left(\frac{x}{2^n}\right) \right\|  \\ && \le \lim_{n\to\infty}2^n  \left\|s\left( 2 g\left(\lambda \frac{x+y}{2^{n+1}}\right)-\lambda g\left(\frac{x}{2^n}\right)-\lambda g\left(\frac{y}{2^n}\right)\right)\right\|\\ && + \lim_{n\to\infty}2^n \left\|t \left(  2h\left(\lambda \frac{x+y}{2^{n+1}}\right)+ 2h\left(\lambda \frac{x-y}{2^{n+1}}\right)- 2\lambda h\left(\frac{x}{2^n}\right)\right) \right\|  +\lim_{n\to\infty} 2^n \varphi\left(\frac{x}{2^n}, \frac{y}{2^n}\right) \\ && = \left\|s\left( 2 D\left(\lambda \frac{x+y}{2}\right)-\lambda D(x)-\lambda D(y)\right)\right\|+ \left\|t \left(  2H\left(\lambda \frac{x+y}{2}\right)+ 2H\left(\lambda \frac{x-y}{2}\right)- 2\lambda H(x)\right) \right\|  
\end{eqnarray*} for  all $\lambda \in {{\mathbb{T}}^1} $ and all $x, y \in A$. So \begin{eqnarray}\label{2.8}&& \|D\left(\lambda(x+y)\right) -\lambda D(x) - \lambda D(y)\| +\| H(\lambda(x+y)) + H(\lambda(x-y) )-2 \lambda H(x) \|   \\ && \le  \left\|s\left( 2 D\left(\lambda \frac{x+y}{2}\right)-\lambda D(x)-\lambda D(y)\right)\right\|+ \left\|t \left(  2H\left(\lambda \frac{x+y}{2}\right)+ 2H\left(\lambda \frac{x-y}{2}\right)- 2\lambda H(x)\right) \right\|  \nonumber  \end{eqnarray} for  all $\lambda \in {{\mathbb{T}}^1} $ and  all $x, y \in A$. 

Let $\lambda =1$ in (\ref{2.8}). 
 By Lemma \ref{lm2.1}, the mappings $D, H: A\to A$ are additive. 

It follows from (\ref{2.8}) and the additivity of $D$ and $H$ that 
\begin{eqnarray*}
 &&\left\|D\left(\lambda (x+y)\right) -\lambda D(x) -\lambda D(y)\right\| + \left\| H\left(\lambda (x+y)\right) - H(\lambda  (x-y)) -2 \lambda H(y) \right\|   \\ && \le \left\| s\left(  D\left(\lambda (x+y)\right) -\lambda D(x) -\lambda D(y)\right)\right\|+ \|t (   H\left(\lambda (x+y)\right) - H(\lambda  (x-y)) -2 \lambda H(y) )\|     \nonumber
\end{eqnarray*} for all $\lambda \in {{\mathbb{T}}^1} $ and all $x, y \in A$.  
Since $|s|< 1$ and $|t|<1$,  \begin{eqnarray*}
D\left(\lambda (x+y)\right) -\lambda D(x) -\lambda D(y) & = & 0,\\
H\left(\lambda (x+y)\right) -  H(\lambda  (x-y)) -2 \lambda H(y) &=& 0
\end{eqnarray*}
and so $D(\lambda x) = \lambda D(x)$ and $H(\lambda x) = \lambda H(x)$  for all $\lambda \in {{\mathbb{T}}^1} $ and all $x, y \in A$.  Thus by Lemma \ref{lm2.2}, the additive  mappings 
$D, H: A\to A$ are $\mathbb{C}$-linear.

It follows from (\ref{2.4}) and the additivity of $D$ and $ H$ that 
\begin{eqnarray*}
&& \|[D,H](xy) - [D,H](x) y - x [D,H](y) \| + \|H(xy) - H(x)y- xH(y)\| 
\\ &&  = 4^n \left\|[g,h]\left(\frac{xy}{4^n}\right) - [g,h]\left(\frac{x}{2^n}\right) \cdot \frac{y}{2^n} - \frac{x}{2^n} \cdot  [g,h]\left(\frac{y}{2^n}\right) \right\|  \\ && + 4^n \left\|h\left(\frac{xy}{4^n}\right) - h\left(\frac{x}{2^n}\right) \cdot \frac{y}{2^n} - \frac{x}{2^n} \cdot h\left(\frac{y}{2^n}\right)\right\|\le 4^n \varphi\left(\frac{x}{2^n},\frac{y}{2^n}\right),
\end{eqnarray*}
which tends to zero as $n\to \infty$, by (\ref{2.2}). So 
\begin{eqnarray*}
[D,H](xy) - [D,H](x) y - x [D,H](y) &=& 0,\\
H(xy) - H(x)y- xH(y) &=& 0
\end{eqnarray*}
for all $x,y\in A$. Hence the mappings  
$[D,H]: A\to A$ and $H : A\to A$ are  derivations.
\end{pf}

\begin{cor}
Let $r > 2$ and $\theta$ be nonnegative real numbers and $g, h :
A \rightarrow A$ be  mappings  satisfying $g(0)=h(0)=0$  and 
\begin{eqnarray}\label{2.9} 
 &&\left\|g\left(\lambda (x+y)\right) -\lambda g(x) -\lambda g(y)\right\| + \left\| h\left(\lambda (x+y)\right)+h(\lambda  (x-y)) - 2\lambda h(x) \right\|  \nonumber  \\ &&\qquad  \le \left\| s\left( 2 g\left(\lambda \frac{x+y}{2}\right)-\lambda g(x)- \lambda g(y)\right)\right\|\\ && \qquad + \left\|t \left(  2h\left(\lambda \frac{x+y}{2}\right)+ 2 h\left(\lambda \frac{x-y}{2}\right)-2\lambda h(x)\right)\right\| +\theta(\|x\|^r + \|y\|^r )    \nonumber
\end{eqnarray} 
for all $\lambda \in {{\mathbb{T}}^1} $ and all $x, y \in A$.  
If the mappings  $g, h: A\to A$ satisfy 
\begin{eqnarray} \label{2.10}
\|[g,h](xy) - [g,h](x) y - x [g,h](y) \| + \|h(xy) - h(x)y- x h(y)\| \le \theta(\|x\|^r + \|y\|^r )   
\end{eqnarray}
 for all $x,y\in A$, then  there exist a unique  $\mathbb{C}$-linear  $D : A \rightarrow A$ and a unique  derivation  $H: A\to A$ such that $[D,H] : A\to A$ is a  derivation and 
\begin{eqnarray*}
\|g(x)-D(x)\| + \|h(x)-H(x)\| \le \frac{2\theta}{2^r -2 } \|x\|^r 
\end{eqnarray*}
for all $x\in A$.
\end{cor}

\begin{pf}
The proof follows from Theorem \ref{thm2.3} by $\varphi (x,y)= \theta (\|x\|^r+\|y\|^r)$  for all $x, y\in A$. 
\end{pf}

\begin{thm}\label{thm2.5}
Let $\varphi : A^2\to [0,\infty)$ be a function  and  $g, h :
A \rightarrow A$ be  mappings satisfying $g(0)=h(0)=0$, {\rm (\ref{2.3})},  {\rm (\ref{2.4})} and 
\begin{eqnarray}\label{2.12}
\Phi(x,y):=  \sum_{j=0}^{\infty} \frac{1}{2^j} \varphi(2^j x, 2^j y) <\infty
\end{eqnarray}
for all $x, y \in A$.
Then there exist a unique  $\mathbb{C}$-linear  $D : A \rightarrow A$ and a unique  derivation  $H: A\to A$ such that $[D,H] : A\to A$ is a  derivation and 
\begin{eqnarray}\label{2.13}
\|g(x)-D(x)\| + \|h(x)-H(x)\| \le \frac{1}{2} \Phi(x,x)
\end{eqnarray}
for all $x\in A$. 
\end{thm}

\begin{pf}
It follows from (\ref{2.6}) that 
\begin{eqnarray}\label{2.14}
\left\|g(x) - \frac{1}{2}g(2x)\right\| + \left\|h(x) - \frac{1}{2}h(2x)\right\|\le \frac{1}{2}\varphi(x,x)
\end{eqnarray}
for all $x\in A$. Thus  
\begin{eqnarray}\label{2.15}
&& \left\|\frac{1}{2^{l} } g\left(\frac{x}{2^l}\right) - \frac{1}{2^{m}} g\left({2^{m}}x\right)\right\| + \left\|\frac{1}{2^{l}} h\left({2^l}x\right) - \frac{1}{2^{m}} h\left({2^{m}}x\right)\right\|    \\ & &  \qquad \le 
 \sum_{j=l}^{m-1}\left\| \frac{1}{2^j } g\left({2^j}x\right) - \frac{1}{2^{j+1}} g\left({2^{j+1}}x\right) \right\| + \sum_{j=l}^{m-1}\left\| \frac{1}{2^j } h\left({2^j}x\right) - \frac{1}{2^{j+1} }h\left({2^{j+1}}x\right) \right\| \nonumber \\ && \qquad 
\le  \nonumber \frac{1}{2} \sum_{j=l}^{m-1} \frac{1}{2^{j}}  \varphi\left({2^{j}}x, {2^{j}}x\right) \
\end{eqnarray}
 for all nonnegative
integers $m$ and $l$ with $m>l$ and all $x \in A$. It follows from
{\rm (\ref{2.15})} that the sequences
$\{\frac{1}{2^{k}} g({2^k}x)\}$ and $\{\frac{1}{2^{k}} h({2^k}x)\}$   are  Cauchy  for all $x \in A$.
Since $Y$ is a  Banach space, the sequences $\{\frac{1}{2^{k} }g({2^k}x)\}$ and $\{\frac{1}{2^{k}} h({2^k}x)\}$   converge.
So one can define
the mappings $D, H : A \rightarrow A$ by
$$D(x) : = \lim_{k\to \infty} \frac{1}{2^{k}} g\left({2^k}x\right),  \qquad H(x) : = \lim_{k\to \infty} \frac{1}{2^{k}} h\left({2^k}x\right)  $$
for all $x \in A$. Moreover, letting $l =0$ and passing to the limit $m \to
\infty$ in {\rm (\ref{2.15})}, we get {\rm (\ref{2.13})}.

By the same reasoning as in the proof of Theorem \ref{thm2.3}, one can show that the mappings   $D, H: A \to A$ are  $\mathbb{C}$-linear.

It follows from (\ref{2.4}) and the additivity of $D, H$ that 
\begin{eqnarray*}
&& \|[D,H](xy) - [D,H](x) y - x [D,H](y) \| + \|H(xy) - H(x)y - xH(y)\| 
\\ &&  = \frac{1}{4^n} \left\|[g,h]\left(4^n xy\right) - [g,h]\left(2^n x\right) (2^n y) - (2^n x) [g,h]\left(2^n y\right) \right\|  \\ && + \frac{1}{4^n} \left\|h\left(4^n xy \right) - h\left(2^n x\right)\cdot 2^n y - 2^n x\cdot h\left(2^n y\right)\right\| \le \frac{1}{4^n} \varphi\left(2^n x, 2^n y\right),
\end{eqnarray*}
which tends to zero as $n\to \infty$, by (\ref{2.12}). So 
\begin{eqnarray*}
[D,H](xy) - [D,H](x) y - x [D,H](y)  &=& 0,\\
H(xy) - H(x)y - xH(y) &=& 0
\end{eqnarray*}
for all $x,y\in A$.  Hence the mappings  
$[D,H]: A\to A$ and $H : A\to A$ are  derivations.
\end{pf}

\begin{cor}
Let $r <1$ and $\theta$ be nonnegative real numbers and $g, h :
A \rightarrow A$ be  mappings  satisfying $g(0)=h(0)=0$,  {\rm (\ref{2.9})} and {\rm (\ref{2.10})}.
Then there exist a unique  $\mathbb{C}$-linear  $D : A \rightarrow A$ and a unique  derivation  $H: A\to A$ such that $[D,H] : A\to A$ is a  derivation and 
\begin{eqnarray*}
\|g(x)-D(x)\| + \|h(x)-H(x)\| \le \frac{2\theta}{2-2^r} \|x\|^r 
\end{eqnarray*}
for all $x\in A$.
\end{cor}

\begin{pf}
The proof follows from Theorem \ref{thm2.5} by $\varphi (x,y)= \theta (\|x\|^r+\|y\|^r)$  for all $x, y\in A$. 
\end{pf}

\section{Stability of additive-additive $(s,t)$-functional inequality {\rm (\ref{0.1})}: a fixed point  method}

Using the fixed point method,  we prove the Hyers-Ulam stability of Lie bracket derivation-derivations   in complex Banach algebras associated to  the additive-additive $(s,t)$-functional
inequality (\ref{0.1}) and the functional inequality (\ref{0.2}).

\begin{thm}\label{thm3.1}
Let $\varphi: A^2 \to [0,\infty)$ be a function  such that 
there exists an $L<1$ with
\begin{eqnarray}\label{3.1}  
 \varphi\left(\frac{x}{2}, \frac{y}{2}\right) 
 \le \frac{L}{4} \varphi\left({x}, y\right)  \le \frac{L}{2} \varphi\left({x}, y\right) 
\end{eqnarray}
 for all $x, y \in A$.
Let $g,h  :
A \rightarrow A$ be  mappings satisfying  $g(0)=h(0)=0$, {\rm (\ref{2.3})} and  {\rm (\ref{2.4})}. Then there exist a unique  $\mathbb{C}$-linear  $D : A \rightarrow A$ and a unique  derivation  $H: A\to A$ such that $[D,H] : A\to A$ is a  derivation and 
\begin{eqnarray}\label{3.2}
\|g(x)-D(x)\| + \|h(x)-H(x)\| \le \frac{L}{2(1-L)} \varphi\left(x,x\right) 
\end{eqnarray}
for all $x\in A$. 
\end{thm}

\begin{pf}
It follows from (\ref{3.1}) that 
$$\sum_{j=1}^{\infty} 4^j \varphi\left( \frac{x}{2^j}, \frac{y}{2^j}\right) \le \sum_{j=1}^{\infty} 4^j \frac{L^j}{4^j} \varphi(x,y) = \frac{L}{1-L} \varphi(x, y) <\infty
$$ for all $x,y \in A$. 
By Theorem \ref{thm2.3}, there exist a unique  $\mathbb{C}$-linear  $D : A \rightarrow A$ and a unique  derivation  $H: A\to A$ satisfying (\ref{2.5}) such that $[D,H] : A\to A$ is a  derivation.

Letting $\lambda =1$ and  $y = x$  in {\rm (\ref{2.3})}, we get
\begin{eqnarray}\label{3.3}
 \| g(2x) - 2g(x) \|+ \|h(2x)-2h(x)\| \le \varphi(x,x)
\end{eqnarray}
for all $x \in A$. 

Consider the set $$ S :  = \{ (g,h) \mid g, h : A \to  A, \ \ g(0)=h(0)=0  \}$$ and introduce the  generalized metric  on $S$:
\begin{eqnarray*}
d((g, h), (g_1, h_1 )) = \inf \left\{ \mu \in {\mathbb R}_+ :\| g(x)-g_1(x)\| + \| h(x)-h_1(x)\|   \le \mu \varphi\left({x}, x\right)  , ~~ \forall x \in A \right\},
\end{eqnarray*} where, as usual, $\inf \phi = +\infty$.
It is easy to show that $(S, d)$ is complete (see  \cite{mr}).

Now we consider the linear mapping $J: S \rightarrow S$ such that
$$J (g,h) (x): = \left(2 g\left(\frac{x}{2}\right) , 2 h\left(\frac{x}{2}\right)\right) $$
for all $x \in A$.

Let $(g, h), (g_1, h_1) \in S$ be given such that $d((g, h), (g_1, h_1)) =  \varepsilon$. Then
$$\|g(x)-g_1(x)\| + \|h(x)-h_1(x)\|\le  \varepsilon \varphi\left({x}, x\right)$$
for all $x\in A$.  Since 
\begin{eqnarray*}
  \left\|2 g\left(\frac{x}{2} \right) - 2 g_1\left( \frac{x}{2}\right)\right\| + \left\|2 h\left(\frac{x}{2} \right) - 2 h_1 \left( \frac{x}{2}\right)\right\|
\le  2\varepsilon  \varphi\left(\frac{x}{2}, \frac{x}{2}\right)  \le 2  \varepsilon \frac{L}{2} \varphi\left({x}, x\right)  = L  \varepsilon \varphi\left({x}, x\right) 
\end{eqnarray*}
for all $x\in A$, 
 $d(J (g,h), J (g_1, h_1))\le L\varepsilon$. This means that
$$d(J (g,h), J(g_1, h_1) ) \le L d((g, h), (g_1, h_1))
$$ for all $(g, h), (g_1, h_1) \in S$.

It follows from {\rm (\ref{3.3})}  that $$\left\| g(x) - 2g\left( \frac{x}{2}\right) \right\|+ \left\| h(x) - 2h\left( \frac{x}{2}\right) \right\| \le \varphi\left( \frac{x}{2}, \frac{x}{2}\right) \le \frac{L}{2} \varphi(x,x)$$ for all $x\in A$. So $d((g,h), (Jg, Jh)) \le \frac{L}{2}$.

By Theorem \ref{thm1.1}, there exist  mappings $D, H : A\rightarrow A$ satisfying the following:

(1) $(D, H)$ is a fixed point of $J$, i.e.,
\begin{eqnarray}\label{3.4}
D\left(x\right) =2D \left(\frac{x}{2}\right), \qquad H\left(x\right) =2H \left(\frac{x}{2}\right)
\end{eqnarray}
 for all $x \in A $.
The mapping $(D,H)$ is a unique fixed point of $J$ in the set $$M=\{g\in
S : d((g, h), (g_1, h_1)) < \infty\} .$$ This implies that $(D,H)$ is a unique mapping
satisfying {\rm (\ref{3.4})} such that there exists a $\mu \in (0, \infty)$
satisfying
\begin{eqnarray*} \|g(x)- D(x)\| + \|h(x)- H(x)\|& \le & \mu   \varphi\left({x}, x\right)
  \end{eqnarray*} for all $x\in A$;

(2) $d(J^l (g,h), (D,H)) \rightarrow 0$ as $l \rightarrow \infty$. This implies the equality

\begin{eqnarray*}
\lim_{l\to \infty} 2^{l}g\left(\frac{x}{2^{l}}\right) = D(x), \quad \lim_{l\to \infty} 2^{l}h\left(\frac{x}{2^{l}}\right) = H(x)
\end{eqnarray*}
for all $x \in A$;

(3) $d((g,h), (D,H)) \le \frac{1}{1-L} d((g,h), J(g,h))$, which implies 
\begin{eqnarray*}
\left\|g(x)- D(x)\right\|+ \left\|h(x)- H(x)\right\|   \le  \frac{L}{2(1-L)}  \varphi\left({x}, x\right) 
\end{eqnarray*}
 for all $x \in A$. 

The rest of the proof is the same as in the proof of Theorem \ref{thm2.3}.
\end{pf}

\begin{cor}
Let $r > 2$ and $\theta$ be nonnegative real numbers and $g, h :
A \rightarrow A$ be  mappings  satisfying $g(0)=h(0)=0$, {\rm (\ref{2.9})} and {\rm (\ref{2.10})}. 
Then there exist  a unique  $\mathbb{C}$-linear  $D : A \rightarrow A$ and a unique  derivation  $H: A\to A$ such that $[D,H] : A\to A$ is a  derivation and 
\begin{eqnarray*}
\|g(x)-D(x)\| + \|h(x)-H(x)\| \le \frac{2\theta}{2^r -2 } \|x\|^r 
\end{eqnarray*}
for all $x\in A$.
\end{cor}

\begin{pf}
The proof follows from Theorem \ref{thm3.1} by taking $L= 2^{1-r}$ and  $\varphi (x,y)= \theta (\|x\|^r+\|y\|^r)$  for all $x, y\in A$. 
\end{pf}

\begin{thm}\label{thm3.3}
Let $\varphi: A^2 \to [0,\infty)$ be a function   such that
there exists an $L<1$ with
\begin{eqnarray}\label{3.5}
 \varphi\left(x, y\right) 
 \le 4 L \varphi\left(\frac{x}{2}, \frac{y}{2}  \right)
\end{eqnarray}
for all $x, y\in A$.  Let $g,h  :
A \rightarrow A$ be  mappings satisfying  $g(0)=h(0)=0$, {\rm (\ref{2.3})} and  {\rm (\ref{2.4})}. Then there exist a unique  $\mathbb{C}$-linear  $D : A \rightarrow A$ and a unique  derivation  $H: A\to A$ such that $[D,H] : A\to A$ is a  derivation and 
\begin{eqnarray}\label{3.6}
\|g(x)-D(x)\| + \|h(x)-H(x)\| \le \frac{1}{2(1-L)} \varphi\left(x,x\right) 
\end{eqnarray}
for all $x\in A$. 
\end{thm}

\begin{pf}
It follows from (\ref{3.5}) that 
$$\sum_{j=1}^{\infty} \frac{1}{4^j} \varphi\left( 2^j x,2^j y\right) \le \sum_{j=1}^{\infty} \frac{1}{4^j} (4 L)^j \varphi(x,y) = \frac{L}{1-L} \varphi(x, y) <\infty
$$ for all $x,y \in A$. 
By Theorem \ref{thm2.5}, there exist a unique  $\mathbb{C}$-linear  $D : A \rightarrow A$ and a unique  derivation  $H: A\to A$   satisfying (\ref{2.13}) such that $[D,H] : A\to A$ is a  derivation.

Let $(S, d)$ be the generalized metric space defined in the proof of Theorem \ref{thm3.1}.

Now we consider the linear mapping $J: S \rightarrow S$ such that
$$J (g,h)(x): = \left( \frac{1}{2} g\left(2x\right) , \frac{1}{2} h\left(2x\right)\right) $$
for all $x \in A$.

It follows from (\ref{3.3}) that 
\begin{eqnarray*}
 \left\| g(x) - \frac{1}{2} g(2x) \right\|+ \left\|h(x)-\frac{1}{2} h(2x)\right\| \le \frac{1}{2} \varphi(x,x)
\end{eqnarray*}
for all $x \in A$.  

The rest of the proof is similar to the proof of Theorem \ref{thm3.1}.
\end{pf}

\begin{cor}
Let $r <1$ and $\theta$ be nonnegative real numbers and $g, h :
A \rightarrow A$ be  mappings  satisfying $g(0)=h(0)=0$,  {\rm (\ref{2.9})} and {\rm (\ref{2.10})}.
Then there exist a unique  $\mathbb{C}$-linear  $D : A \rightarrow A$ and a unique  derivation  $H: A\to A$ such that $[D,H] : A\to A$ is a  derivation and 
\begin{eqnarray*}
\|g(x)-D(x)\| + \|h(x)-H(x)\| \le \frac{2\theta}{2-2^r} \|x\|^r 
\end{eqnarray*}
for all $x\in A$.
\end{cor}

\begin{pf}
The proof follows from Theorem \ref{thm3.3} by taking $L= 2^{r-1}$ and  $\varphi (x,y)= \theta (\|x\|^r+\|y\|^r)$  for all $x, y\in A$. 
\end{pf}

\section*{Acknowledgments}

 C. Park  was supported by Basic Science Research Program through the National Research Foundation of Korea funded by the Ministry of Education, Science and Technology (NRF-2017R1D1A1B04032937).

\section*{Competing interests}

The author declares that they have no competing interests.

\section*{Authors' contributions}

The author conceived of the study, participated in its design and
coordination, drafted the manuscript, participated in the sequence
alignment, and read and approved the final manuscript.

{\footnotesize


\begin{thebibliography}{99}


 \bibitem {a} T. Aoki,
\newblock{\it On the stability of the linear transformation in Banach spaces},
 \newblock J. Math. Soc. Japan {\bf 2} (1950), 64--66.




\bibitem {cr1}  L. C\u{a}dariu, V.   Radu, {\it  Fixed points and the stability of Jensen's functional equation},
\newblock J. Inequal. Pure Appl. Math. {\bf 4}, no. 1, Art. ID 4 (2003).

\bibitem {cr2}  L. C\u{a}dariu, V.  Radu, {\it   On the stability of the Cauchy functional equation: a fixed point approach}, 
 \newblock Grazer Math. Ber. {\bf 346}  (2004), 43--52.

 \bibitem {cr3} L.  C\u{a}dariu, V.  Radu, {\it   Fixed point methods for the generalized stability of functional equations in a single variable}, 
 \newblock Fixed Point Theory  Appl. {\bf 2008}, Art. ID 749392 (2008).



\bibitem {dm}  J. Diaz, B.  Margolis, {\it   A fixed point theorem of the alternative for contractions on a generalized
 complete metric space},   \newblock Bull. Am. Math. Soc. {\bf 74} (1968), 305--309.



\bibitem{e17} I. EL-Fassi,  {\it Generalized hyperstability of a Drygas functional equation on a restricted domain using Brzdek's fixed point theorem,} J.  Fixed Point Theory Appl.  {\bf 19} (2017), 2529--2540.

  \bibitem {el} I. EL-Fassi, {\it   Solution and approximation of radical quintic functional equation related to quintic mapping in quasi-$\beta$-Banach spaces}, 
Rev. R. Acad. Cienc. Exactas F\'{i}s. Nat. Ser. A Mat.  {\bf 113} (2019), no. 2, 675--687.



\bibitem{ega} M. Eshaghi Gordji, M.B. Ghaemi and B. Alizadeh,  {\it A fixed point method for perturbation of higher ring derivationsin non-Archimedean Banach algebras,} Int. J. Geom.
 Meth. Mod. Phys. {\bf 8} (2011), no. 7,   1611--1625.



\bibitem{ga94} P. G\v avruta,
\newblock{\it A generalization of the Hyers-Ulam-Rassias stability of
approximately additive mappings},
\newblock  J. Math. Anal. Appl. {\bf 184} (1994), 431--436.



\bibitem {hy41} D.H. Hyers,
\newblock{\it On  the stability of the linear functional equation},
 \newblock Proc. Nat. Acad. Sci. U.S.A. {\bf 27} (1941), 222--224.

  \bibitem {ir} G.  Isac, Th. M. Rassias, {\it   Stability of $\psi$-additive mappings: Applications to nonlinear analysis}, 
Int. J. Math. Math. Sci. {\bf 19} (1996), 219--228.


 \bibitem {mr}  D. Mihe\c{t}, V.  Radu, {\it   On the stability of the additive Cauchy functional equation in random normed spaces}, 
 J. Math. Anal. Appl. {\bf 343} (2008), 567--572.

 
\bibitem {n15} I. Nikoufar, {\it Jordan $(\theta, \phi)$-derivations on Hilbert $C^*$-modules}, 
Indag. Math. {\bf 26} (2015), 421--430.




\bibitem {pa06} C. Park, \newblock {\it Homomorphisms between Poisson
$JC^*$-algebras}, \newblock Bull. Braz. Math. Soc. {\bf 36}
(2005), 79--97.


\bibitem {ppp} C. Park, {\it Additive $\rho$-functional inequalities and equations}, J. Math. Inequal. {\bf 9} (2015), 17--26.

\bibitem {ppp1} C. Park, {\it Additive $\rho$-functional inequalities in non-Archimedean normed spaces}, J. Math. Inequal. {\bf 9} (2015), 397--407.



\bibitem {pa17} C. Park, {\it Fixed point method for set-valued functional equations}, J. Fixed Point Theory Appl.  {\bf 19} (2017),  2297--2308.




\bibitem {par} C. Park, {\it Biderivations and bihomomorphisms in Banach algebras}, Filomat (in press).



  \bibitem {r}  V. Radu, {\it The fixed point alternative and the stability of functional equations},  \newblock
Fixed Point Theory  {\bf 4} (2003), 91--96.

 \bibitem {ra78} Th. M. Rassias,
\newblock{\it  On  the stability of the linear mapping in Banach spaces},
 \newblock Proc. Am. Math. Soc. {\bf 72} (1978), 297--300.

  \bibitem {sl}L. Sz\'{e}kelyhidi, {\it   Superstability of functional equations related to spherical functions}, 
Open Math.  {\bf 15} (2017), no. 1, 427--432.

\bibitem {ul60} S. M. Ulam,
 \newblock{\it A Collection of the Mathematical Problems},
 \newblock  Interscience Publ. New York, 1960.


\bibitem {w} Z. Wang,  {\it Stability of two types of cubic fuzzy set-valued functional equations,}  Results Math. {\bf 70} (2016), 1--14.



\end{thebibliography}
\end{document}